\documentclass[twoside,12pt]{article}
\usepackage[all, 2cell, dvips]{xy}
\usepackage{latexsym,amsmath, amsfonts, graphics, epsf, epic}
\usepackage{amsthm}
\usepackage{amssymb}
\usepackage{amsopn}
\usepackage{amscd}
\usepackage{color}

\theoremstyle{plain}
\newtheorem{thm}{Theorem}[section]

\newtheorem{lem}[thm]{Lemma}
\newtheorem{prop}[thm]{Proposition}

\theoremstyle{definition}

\newtheorem{remark}[thm]{Remark}

\usepackage{amssymb}

\def\lm{\lambda}

\def\l.l.o.{\it l.l.o}

\def\chiup{\raise 2pt\hbox{$\chi$}}

\title{Humps for Dyck and for Motzkin paths}

\author{ A. Regev\\Department of Mathematics,\\ The Weizmann Institute of Science,
\\Rehovot, Israel\\e-mail: amitai.regev~at~weizmann.ac.il}

\date{}

\begin{document}

\maketitle

\begin{abstract}
   We calculate the total number of humps in Dyck and in Motzkin
paths, and we give Standard-Young-Tableaux-interpretations of the
numbers involved. One then observes the intriguing phenomena that
the humps-calculations change the partitions in a strip to
partitions in a hook.

\end{abstract}

\section{Introduction}

Let $\lm$ be a partition and denote by $f^\lm$ the number of
standard Young tableaux (SYT) of shape $\lm$. The number $f^\lm$
can be computed  for example by the hook formula~\cite[Corollary
7.21.6]{stanley}.

\subsection{The $(k,\ell)$-hook sums}
We consider the SYT in the $(k,\ell)$ hook. More precisely, given
integers $k,\ell,n\ge 0$ we denote
\[
H(k,\ell;n)=\{\lm=(\lm_1,\lm_2,\ldots)\mid \lm\vdash n~\mbox{and}~
\lm_{k+1}\le \ell\}\qquad\mbox{and}\qquad S(k,\ell;n)=\sum_{\lm\in
H(k,\ell;n)}f^\lm.
\]
We remark that classically, the partitions $\lm\in \cup_{n\ge
0}H(k,0;n)$ parametrize the irreducible representations of the Lie
algebra $gl(k,\mathbb{C})$. Also, the partitions $\lm\in
\cup_{n\ge 0}H(k,\ell;n)$ parametrize those  of the Lie
super-algebra $pl(k,\ell)$~\cite{berele}.


\medskip

For the "strip" sums $S(k,0;n)$ it is
known~\cite{R1}~\cite{stanley} that
\[
S(2,0;n)={n\choose\lfloor\frac{n}{2}\rfloor}\quad\mbox{and}\quad
S(3,0;n)=\sum_{j\ge 0}\frac{1}{j+1}{n\choose 2j}{2j\choose j}.
\]
 Furthermore,
Gouyon-Beauchamps~\cite{gouyon}~\cite{stanley} proved that
\[
S(4,0;n)=C_{\lfloor\frac{n+1}{2}\rfloor}\cdot
C_{\lceil\frac{n+1}{2}\rceil}\quad\mbox{and}\quad
S(5,0;n)=6\sum_{j=0}^{\lfloor\frac{n}{2}\rfloor}{n\choose 2j}\cdot
C_j\cdot\frac{(2j+2)!}{(j+2)!(j+3)!},
\]
where
\begin{eqnarray}\label{catalan1}
C_j=\frac{1}{j+1}{2j\choose j}\qquad\mbox{are the Catalan
numbers.}
\end{eqnarray}

\medskip
So far only the "hook" sums $S(1,1;n)$ and $S(2,1;n)=S(1,2;n)$
have been calculated: it is easy to see that $S(1,1;n)=2^{n-1}$;
for the sum $S(2,1;n)$ see Equation~\eqref{motzkin.path.3} below.

\subsection{Humps-calculations of paths}

\subsubsection{Dyck paths}

The Catalan numbers $C_n$ count many combinatorial
objects~\cite{stanley}. For example, they count the number of SYT
of the $2\times n$ rectangular shape $(n,n)$, namely
$C_n=f^{(n,n)}$, which can be thought of as a SYT-interpretation
of the Catalan numbers. It is also well known that $C_n$ counts
the number of Dyck paths of length $2n$. A Dyck path of length
$2n$ is a lattice path, in $\mathbb{Z}\times \mathbb{Z}$, from
$(0,0)$ to $(2n,0)$, using up-steps $(1,1)$ and down-steps
$(1,-1)$ and never going below the $x$-axis.

\medskip
 A {\it hump} in a Dyck path is an up-step followed by  a
down-step. For example, there are 2 Dyck paths of length 4, with
total number of humps being 3; and there are 5 Dyck paths of
length 6, with total number of humps being 10. We denote by ${\cal
H}C_n$ the total number of humps in the Dyck paths of length $2n$.
Thus ${\cal H}C_2=3$ and ${\cal H}C_3=10$.
Remark~\ref{dyck.humps.11}.1 gives a SYT-interpretation of the
numbers ${\cal H}C_n$.

\medskip

\subsubsection{Motzkin paths}
A  Motzkin path of length
$n$ is a lattice path from $(0,0)$ to $(n,0)$, using flat-steps
$(1,0)$, up-steps $(1,1)$ and down-steps $(1,-1)$, and never going
below the $x$-axis. The Motzkin number $M_n$ counts the number of
Motzkin paths of length $n$. Recall also~\cite{R1}~\cite{stanley}
that
\begin{eqnarray}\label{motzkin.path.1}
M_n=S(3,0;n)=\sum_{\lm\in H(3,0;n)}f^\lm,
\end{eqnarray}
which gives a SYT-interpretation of the Motzkin numbers $M_n$.
Bijective proofs of~\eqref{motzkin.path.1}, as well as of related
identities due to Zeilberger~\cite{doron}, have recently been
given in~\cite{eu}.

\medskip
A hump in a Motzkin path is an up-step followed by zero or more
flat-steps followed by a down-step. We denote by ${\cal H}M_n$ the
total number of humps in the Motzkin paths of length $n$. Thus
${\cal H}M_2=1$ and ${\cal H}M_3=3$. Theorem~\ref{motzkin.humps.1}
gives a SYT-interpretation of the numbers ${\cal H}M_n$.

\subsubsection{The main results}
The main results here are explicit formulas for the total number
of humps for the Dyck paths and for the Motzkin paths of a given
length, together with SYT-interpretation of these numbers.
 One then observes the following intriguing
phenomena. Theorem~\ref{dyck.humps.1} below shows that
 the total number
of humps ${\cal H}C_n$ for the Dyck paths of length $2n$ satisfies
${\cal H}C_n=f^{(n,1^{n})}$. Together with $C_n=f^{(n,n)}$ this
shows, roughly, that the humps-calculations correspond the
$2\times n$ rectangular shape $(n,n)$ to the $1-1$ hook shape
$(n,1^n)$. A somewhat similar phenomena occurs when studying humps
in Motzkin paths: while $M_n=S(3,0;n)$,
Theorem~\ref{motzkin.humps.1} asserts that $ {\cal
H}M_n=S(2,1;n)-1,$ which gives a SYT-interpretation of the numbers
${\cal H}M_n$. This shows, roughly, that the humps-calculations
correspond the $(3,0)$ strip shape to the $(2,1)$  hook shape.

\medskip

We also consider "super" such paths, which are allowed to also go
below the $x$-axis, and  show that their number is essentially
twice the number of the corresponding humps, see
Remarks~\ref{dyck.humps.11}.2 and~\ref{motzkin.humps.111}.

\medskip
In Section~\ref{double1} we give a double Dyck path interpretation
for the sums $S(4,0;n)$, together with the corresponding
hump-numbers ${\cal H}S(4,0;n)$. Proposition~\ref{double2} gives
the intriguing identity ${\cal H}S(4,0;n)=\frac{n+3}{2}\cdot
S(4,0;n)$.
\begin{remark}
It would be interesting to find bijective proofs to the identities
in this paper.
\end{remark}

\section{Humps for Dyck paths}

Recall that $ C_n=f^{(n,n)}=\frac{1}{n+1}{2n\choose n}, $ which is
a SYT-interpretation of the Catalan numbers $C_n$.
\begin{thm}\label{dyck.humps.1}
Let ${\cal H}C_n$ denote the total number of humps for all the
Dyck paths from $(0,0)$ to $(2n,0)$, then
\begin{eqnarray}\label{late.mozkin.1}
{\cal H}C_n={2n-1\choose n}.
\end{eqnarray}
\end{thm}
\begin{remark}\label{dyck.humps.11}
\begin{enumerate}
\item
Let $\lm$ be the $1-1$ hook shaped diagram $\lm=(n,1^{n})$, then
\[
f^\lm={2n-1\choose n}={\cal H}C_n,
\]
which gives a SYT-interpretation of the numbers ${\cal H}C_n$.
\item
A super Dyck path is a Dyck path which is allowed to also go below
the $x$-axis. Let $SD_n$ denote the number of super Dyck paths of
length $2n$. By standard arguments it follows that
\[
SD_n={2n\choose n}=2\cdot {2n-1\choose n}=2\cdot{\cal H}C_n.
\]
\end{enumerate}

\end{remark}

The proof of Theorem~\ref{dyck.humps.1} clearly follows from the
following two lemmas.

\begin{lem}\label{dyck.humps.2}
We have ${\cal H}C_0={\cal H}C_1=C_0=C_1=1,$ and the Catalan
numbers $C_n$ and the numbers ${\cal H}C_n$ satisfy the following
equation
\begin{eqnarray}\label{dyck.humps.3}
{\cal H}C_n={\cal H}C_{n-1}+\sum_{j=1}^{n-1} ({\cal H}C_{j-1}
\cdot C_{n-j}+
 C_{j-1} \cdot {\cal H}C_{n-j}).
\end{eqnarray}
\end{lem}
\begin{proof}
Given a Dyck path ${\cal D}={\cal D}_n$ of length $2n$, we read it
from left to right. Then ${\cal D}_n$ meets the $x$-axis for the
{\it first time} after $2j$ steps, $1\le j\le n$.

\medskip
{\it Case 1:} $j=n$. This implies that except for the endpoints
$(0,0)$ and $(2n,0)$, ${\cal D}$ does not meet the $x$-axis. Thus
${\cal D}_n=s_1\cdots s_n$ where $s_1$ is an up step, $s_n$ is a
down step, and $s_2\cdots s_{n-1}$ corresponds to a Dyck path
${\cal D}_{n-1}$ in an obvious way. Also, the number of humps of
 $s_1\cdots s_n$ and of $s_2\cdots s_{n-1}$ is the same.

\medskip
It follows that the Dyck paths of Case 1 contribute ${\cal
H}C_{n-1}$ to the total number of haumps ${\cal H}C_{n-1}$.

\medskip
{\it Case 2:} $j\le n-1$. Thus ${\cal D}_n$ is the concatenation
of two Dyck paths ${\cal D}_n={\cal D}_{2j}{\cal D}_{2(n-j)}$,
where ${\cal D}_{2(n-j)}$ is an arbitrary Dyck path of length
$2(n-j)$ but ${\cal D}_{2j}$, of length $2j$, is of the type
studied in Case 1. Thus there are $C_{j-1}$ paths ${\cal D}_{2j}$,
with total humps-contribution being ${\cal H}C_{j-1}$. And there
are $C_{2(n-j)}$ paths ${\cal D}_{2(n-j)}$, with total
humps-contribution being ${\cal H}C_{2(n-j)}$.

\medskip
Equation~\eqref{dyck.humps.3} now follows, since the number of
humps ${\cal H}({\cal D}_n)$  of ${\cal D}_n={\cal D}_{2j}{\cal
D}_{2(n-j)}$ is the sum ${\cal H}({\cal D}_n)={\cal H}({\cal
D}_{2j})+{\cal H}({\cal D}_{2(n-j)}).$

\end{proof}

\begin{lem}\label{dyck.humps.4}
Together with ${\cal H}C_0={\cal H}C_1=1,$
Equation~\eqref{dyck.humps.3} implies
Equation~\eqref{late.mozkin.1}.
\end{lem}

\begin{proof}

[D. Zeilberger] By induction on $n$, replace each ${\cal H}C_k$ on
the right hand side of Equation~\eqref{dyck.humps.3} by
${2k-1\choose k}$. We obtain the following binomial identity:
\begin{eqnarray}\label{dyck.humps.31}{2n-1\choose
n}={2n-3\choose n-1}+\sum_{j=1}^{n-1} \left({2j-3\choose j-1}
\cdot C_{n-j} +C_j
 \cdot {2n-2j-1\choose n-j}\right)
\end{eqnarray}
This identity is easily provavle by Gosper's algorithm for
indefinite summation (implemented in Maple's sum command).

\medskip
For a direct proof note first that ${2k-1\choose
k}=\frac{k+1}{2}\cdot C_k$ if $k\ge 1$. Thus the right side
of~\eqref{dyck.humps.31} is
\[
2{2(n-1)-1\choose n-1}+C_{n-1}+
\sum_{j=1}^{n-1}\left(\frac{j}{2}\cdot C_{j-1}\cdot
C_{n-j}+C_{j-1}\cdot \frac{n-j+1}{2} \cdot
C_{n-j}\right)=~~~~~~~~~~~
\]
\[
~~~~~~~~~~~~~~~~~~~~~~~~~~~=\frac{n+1}{2}\cdot
\sum_{j=1}^nC_{j-1}\cdot C_{n-j} = \frac{n+1}{2}\cdot C(n)
\]
(by the defining relation for $C_n$: $C_n=\sum_{j=1}^nC_{j-1}\cdot
C_{n-j}$), and this equals the left side of~\eqref{dyck.humps.31}.

\end{proof}

This completes the proof of Theorem~\ref{dyck.humps.1}.

\section{Humps for Motzkin paths}

Recall that a  Motzkin path of length $n$ is a lattice path from
$(0,0)$ to $(n,0)$, using flat-steps $(1,0)$, up-steps $(1,1)$ and
down-steps $(1,-1)$ and never going below the $x$-axis. The
Motzkin number $M_n$ counts the number of Motzkin paths of length
$n$. Recall also~\cite{R1}~\cite{stanley} that
\begin{eqnarray}\label{motzkin.path.01}
M_n=\sum_{\lm\in H(3,0;n)}f^\lm,
\end{eqnarray}
 which gives a SYT-interpretation
of the Motzkin numbers $M_n$. A bijective proof
of~\eqref{motzkin.path.1} has recently been given in~\cite{eu}.

\medskip
Recall that a hump in a Motzkin path is an up-step followed by 0
or more flat-steps followed by a down-step.

\begin{thm}\label{motzkin.humps.2}
The Number ${\cal H}M_n$ of humps in all Motzkin paths of length n
is given by
\begin{eqnarray}\label{motzkin.path.2}
{\cal H}M_n=\frac{1}{2}\sum_{j\ge 1}{n\choose j}{n-j\choose j}.
\end{eqnarray}
\end{thm}

\begin{remark}\label{motzkin.humps.111}

A super Motzkin path is a Motzkin path which is allowed to also go
below the $x$-axis. Let $SM_n$ denote the number of super Motzkin
paths of length $n$, then one can prove the recurrence
\[
SM_n=SM_{n-1}+2\cdot\sum_{k=2}^n M_{k-2}\cdot SM_{n-k},
\]
which then implies that
\[
SM_n=\sum_{j\ge 0}{n\choose j} {n-j\choose j}.
\]
Together with Theorem~\ref{motzkin.humps.2} it implies that
\[
SM_n=
2\cdot {\cal H}M_n
+1.
\]
\end{remark}

{\bf Proof of Theorem~\ref{motzkin.humps.2}}
\begin{proof}
This result is stated in~\cite[sequence A097861]{sloane}, and it
can be proved as follows. First argue as in
Lemma~\ref{dyck.humps.2}: Given a Motzkin path ${\cal M}={\cal
M}_n$ of length $n$, we read it from left to right. Then ${\cal
M}_n$ meets the $x$-axis for the {\it first time} after $j$ steps,
$1\le j\le n$. The case $j=1$ contributes ${\cal H}M_{n-1}$ to
${\cal H}M_{n}$. If $2\le j$, then ${\cal M}_n={\cal M}_j{\cal
M}_{n-j}$ (concatenation) where ${\cal M}_j$ starts with an
up-step and ends with a down-step, while ${\cal M}_{n-j}$ is an
arbitrary Motzkin path of length $n-j$. Thus ${\cal M}_j$
corresponds to a Motzkin path ${\cal M}'_{j-2}$ of length $j-2$,
with the same number of humps as ${\cal M}_j$, except for the case
that ${\cal M}'_{j-2}$ is a sequence of flat steps, in which case
${\cal M}_j$ contributes one hump while ${\cal M}'_{j-2}$
contributes zero humps.

\medskip
Now ${\cal H}M_0={\cal H}M_1=0$, and for $n\ge 2$, similar to the
proof of Lemma~\ref{dyck.humps.2}, the above argument implies the
recurrence
\begin{eqnarray}\label{reccurence.motzkin.1}
{\cal H}M_n={\cal H}M_{n-1}+\sum_{k=2}^n\left((1+{\cal
H}M_{k-2})\cdot M_{n-k}+M_{k-2}\cdot {\cal H}M_{n-k}\right).
\end{eqnarray}
Denote
\begin{eqnarray}\label{motzkin.path.22}
B_n=\frac{1}{2}\sum_{j\ge 1}{n\choose j}{n-j\choose j},
\end{eqnarray}
then the proof of Theorem~\ref{motzkin.humps.2} follows -- by
induction on $n$ -- from Equation~\eqref{reccurence.motzkin.1} and
from the following binomial identity
\begin{eqnarray}\label{reccurence.motzkin.011}
B_n=B_{n-1}+\sum_{k=2}^n\left((1+B_{k-2})\cdot
M_{n-k}+M_{k-2}\cdot B_{n-k}\right).
\end{eqnarray}
 Equation~\eqref{reccurence.motzkin.011} can be proved by the WZ
method~\cite{doron3},~\cite{doron2}.
\end{proof}

 Recall that $S(2,1;n)=
\sum_{\lm\in H(2,1:n)}f^\lm.$ Here we prove the following
intriguing identity.

\begin{thm}\label{motzkin.humps.1}
Let ${\cal H}M_n$ denote total number of humps for all the Motzkin
paths of length $n$, then
\[
{\cal H}M_n=S(2,1;n)-1.
\]
\end{thm}
This gives a SYT-interpretation of the numbers ${\cal H}M_n$.

\begin{proof}
 We prove Theorem~\ref{motzkin.humps.1} numerically. The following
 equation
was proved in~\cite[Theorem 8.1]{R2}:
\begin{eqnarray}\label{motzkin.path.3}
S(2,1;n)=~~~~~~~~~~~~~~~~~~~~~~~~~~~~~~~~~~~~~~~~~~~~~~~~~~~~~~~~~~~~~~~~~~~~~~~~~~~~~~~~~~~~~~~~~~~~
\end{eqnarray}
\[
~~~~~~=\frac{1}{4}\left(\sum_{r=0}^{n-1}{n-r\choose{\lfloor\frac{n-r}{2}\rfloor}}
{n\choose r}
+\sum_{k=1}^{\lfloor\frac{n}{2}\rfloor-1}\frac{n!}{k!\cdot
(k+1)!\cdot (n-2k-2)!\cdot (n-k-1)\cdot(n-k)}\right)+1
\]
Combining equations~\eqref{motzkin.path.2}
and~\eqref{motzkin.path.3}, the proof of
Theorem~\ref{motzkin.humps.1} will follow once  the following
binomial identity -- of interest on its own -- is proved.
\begin{lem}\label{motzkin.humps.11}
For $n\ge 2$
\begin{eqnarray}\label{motzkin.path.222}
2\sum_{j\ge 1}{n\choose j}{n-j\choose
j}=~~~~~~~~~~~~~~~~~~~~~~~~~~~~~~~~~~~~~~~~~~~~~~~~~~~~~~~~~~~~~~~~~~~~~~~~~~~~~~~~
\end{eqnarray}
\[
=\sum_{r=0}^{n-1}{n-r\choose{\lfloor\frac{n-r}{2}\rfloor}}
{n\choose r}
+\sum_{k=1}^{\lfloor\frac{n}{2}\rfloor-1}\frac{n!}{k!\cdot
(k+1)!\cdot (n-2k-2)!\cdot (n-k-1)\cdot(n-k)}.
\]
\end{lem}

Equation~\eqref{motzkin.path.222} can be proved by the WZ
method~\cite{doron3},~\cite{doron2}.

\medskip
This completes the proof of Theorem~\ref{motzkin.humps.1}.
\end{proof}

\section{Humps for $S(4,0;n)$}\label{double1}
\subsection{Double Dyck paths interpretation for $S(4,0;n)$}

By~\cite{gouyon},
\begin{eqnarray}\label{dyck.humps.41}
S(4,0;n)=C_{\lfloor\frac{n+1}{2}\rfloor}\cdot
C_{\lceil\frac{n+1}{2}\rceil}.
\end{eqnarray}
Therefore the Dyck paths interpretation for $C_m$ implies the
following interpretation for $S(4,0;n)$:

\medskip
$S(4,0;n)$ is the number of Dyck paths from $(0,0)$ to
$(2(n+1),0)$ which go through $(2\lfloor\frac{n+1}{2}\rfloor,0)$.
Such a path is the concatenation of two Dyck paths, one of length
$2\lfloor\frac{n+1}{2}\rfloor$ and one of length
$2\lceil\frac{n+1}{2}\rceil$   (note that for any integer $m$,
$\lfloor\frac{m}{2}\rfloor+\lceil\frac{m}{2}\rceil=m)$. We call it
a {\it double-Dyck} path. Clearly, the number of humps of such a
concatenated path is the sum of the humps of its two parts. By
arguments similar to those in the proof of
Lemma~\ref{dyck.humps.2} and of Theorem~\ref{motzkin.humps.2} this
implies

\begin{prop}\label{dyck.humps.40}
Let
 ${\cal H}S(4,0;n)$ denote the total number of humps of the double-Dyck paths
 corresponding to $S(4,0;n)$, then
\[
{\cal H}S(4,0;n)={\cal H} C_{\lfloor\frac{n+1}{2}\rfloor}\cdot
C_{\lceil \frac{n+1}{2}\rceil}+
C_{\lfloor\frac{n+1}{2}\rfloor}\cdot {\cal H} C_{\lceil
\frac{n+1}{2}\rceil}.
\]
\end{prop}
For $n=1,2,\ldots$ this gives the following values for ${\cal
H}S(4,0;n)$:

2, 5, 12, 35, 100, 315, 980, 3234, 10584, 36036, 121968, 424710,
etc.

\medskip
We have the following intriguing identity:
\begin{prop}\label{double2}
${\cal H}S(4,0;n)=\frac{n+3}{2}\cdot S(4,0;n)$.
\end{prop}

\begin{proof}
The proof is a straight forward calculations, applying
Theorem~\ref{dyck.humps.1}, Equation~\eqref{dyck.humps.41} and
Proposition~\ref{dyck.humps.40}.
\end{proof}

{\bf Acknowledgment}. We thank D. Zeilberger for some very useful
suggestions.

\end{document}